\documentclass[10pt,reqno]{amsart}
\usepackage{hyperref}
\usepackage{float}
\usepackage{mathtools}
\usepackage{amscd,amssymb,amsmath,latexsym,bm}
\usepackage[mathcal,mathscr]{euscript}
\usepackage{lipsum}
\usepackage{amsfonts}
\usepackage{graphicx}
\usepackage{epstopdf}
\usepackage{algorithmic}
\usepackage{hyperref}
\usepackage{xcolor}
\usepackage{upgreek}
\usepackage{enumerate}
\usepackage{ulem}
\usepackage{stmaryrd}
\usepackage{gensymb}			
\usepackage[utf8]{inputenc}
\usepackage[small,nohug]{diagrams}
\diagramstyle[labelstyle=\scriptstyle]
 
\newtheorem{theorem}{Theorem}[section]

\newtheorem{corollary}[theorem]{Corollary}
\newtheorem{proposition}[theorem]{Proposition}

\newtheorem{definition}[theorem]{Definition}

\numberwithin{equation}{section}
\numberwithin{figure}{section}

\newcommand{\CM}{{\mathbb C}}
\newcommand{\NM}{{\mathbb N}}

\newcommand{\RM}{{\mathbb R}}
\newcommand{\SM}{{\mathbb S}}

\newcommand{\ZM}{{\mathbb Z}}

\newcommand{\KM}{{\mathbb K}}

\newcommand{\Gg}{{\mathcal G}}

\newcommand{\Cc}{{\mathcal C}}

\newcommand{\Ll}{{\mathcal L}}

\newcommand{\Kk}{{\mathcal K}}

\begin{document}

\title{Topological Lattice Defects by Groupoid Methods and Kasparov's $KK$-Theory}

\author{Emil Prodan}

\address{Department of Physics and
\\ Department of Mathematical Sciences 
\\Yeshiva University 
\\New York, NY 10016, USA \\
\href{mailto:prodan@yu.edu}{prodan@yu.edu}}

\date{\today}

\begin{abstract} 
The bulk-boundary and a new bulk-defect correspondence principles are formulated using groupoid algebras. The new strategy relies on the observation that the groupoids of lattices with boundaries or defects display spaces of units with invariant accumulation manifolds, hence they can be naturally split into disjoint unions of open and closed invariant sub-sets. This leads to standard exact sequences of groupoid $C^\ast$-algebras that can be used to associate a Kasparov element to a lattice defect and to formulate an extremely general bulk-defect correspondence principle. As an application, we establish a correspondence between topological defects of a 2-dimensional square lattice and Kasparov's group $KK^1 (C^\ast(\ZM^3),\CM)$. Numerical examples of non-trivial bulk-defect correspondences are supplied.
\end{abstract}

\thanks{This work was supported by the National Science Foundation through the grant DMR-1823800.}

\maketitle


\setcounter{tocdepth}{1}

\section{Introduction and Main Statement}

The framework of $C^\ast$-algebras, with its standard $K-$ and $KK-$theoretic tools, proved to be extremely fruitful for studying the bulk-boundary correspondence principle for topological insulators and super-conductors \cite{KellendonkRMP2002,BourneMPL2015,ProdanSpringer2016,
BourneAHP2017,BourneAHP2020}. Disordered and quasi-periodic lattice systems can be conveniently analyzed in the framework of crossed-product algebras and the bulk-boundary correspondence principle can be derived from the exact sequence of boundary, half-space and bulk algebras, which in essence is just the Pimsner-Voiculescu exact sequence \cite{PimsnerJOP1980}. Roe algebras \cite{HigsonBook} have been also successfully used to explore the bulk-boundary correspondence principle in settings with irregular boundaries \cite{ThiangJGP2020,LudewigArxiv2020} or in hyperbolic geometries \cite{LudwigArxiv2020}. Recently, Roe algebras have been also employed to formulate a bulk-defect correspondence principle for a weak topological insulator \cite{KubotaJPA2021}. However, our focus here will be on groupoid algebras. In many respects, they supply the most natural framework for analyzing discrete resonating systems. Indeed, Bellissard and Kellendonk showed \cite{Bellissard1986,Bellissard1995,KellendonkRMP95} that each uniformly discrete pattern in the euclidean space can be associated with a specific \'etale groupoid $\Gg$ and groupoid $C^\ast$-algebra $C^\ast_r(\Gg)$. Any Galilean invariant time evolution operator over such pattern can be realized from a left-regular representation of this algebra. The framework of groupoid algebras has been successfully employed in the study of topological amorphous systems \cite{BourneJPA2018} and to formulate \cite{BourneAHP2020} a bulk-boundary correspondence principle in a very generic setting.

\begin{figure}[t]
\center
\includegraphics[width=0.5\textwidth]{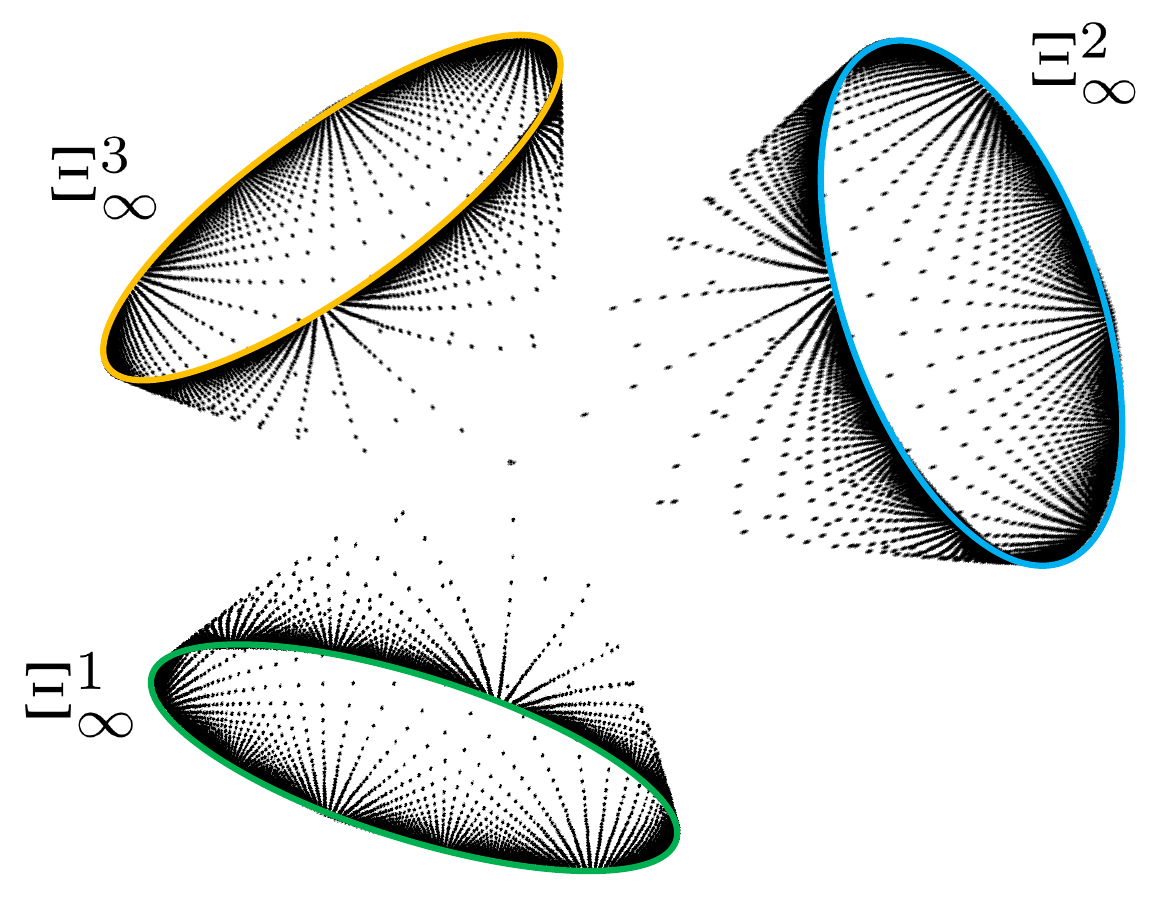}\\
  \caption{\small Example of a space of units occuring in the present work.
  }
 \label{Fig:GenTrans}
\end{figure}

In this work, we present a conceptually different approach to the bulk-boundary principle and formulate a new bulk-defect correspondence principle. To explain the strategy, let us recall that any uniformly discrete pattern leads to a closed topological space $\Xi$ inside the metric space of closed subsets of the euclidean space. This space $\Xi$ serves as the space of units for the groupoid $\Gg$ mentioned above. Traditionally, the bulk-boundary correspondence principle was initiated from the bulk, where one is dealing with a Delone set, not just a uniformly discrete one. This, however limits the array of lattice defects that can be generated, specifically, to boundaries and domain-walls. For example, the lattice defects introduced in Sec.~\ref{Sec:NewExample} cannot be generated with a fixed bulk system. To eliminate this shortcoming, we focus from the beginning on the groupoid associated with the lattice with either a boundary or defect. In the examples worked out here, the space of units $\Xi$ of this groupoid looks like in Fig.~\ref{Fig:GenTrans}. It has closed invariant accumulation manifolds $\Xi_\infty^i$, which result from exploring the physical pattern in different asymptotic regions, infinitely far from the defect and where the pattern looks like a rotated bulk sample. As such, the space of units can be naturally divided as $\Xi=\Xi_\infty^c \cup \Xi_\infty$ ($\Xi_\infty^c$ being the complement of $\Xi_\infty$), or other combinations, which sets an exact sequence of $C^\ast$-algebras \cite{WilliamsBook}
\begin{equation}\label{Eq:ExactSeq0}
\begin{diagram}
0 & \rTo & C_r^\ast(\Gg|_{\Xi_\infty^c}) & \rTo{\rm i \ } & C_r^\ast(\Gg) & \rTo{\rm j \ } & C_r^\ast(\Gg|_{\Xi_\infty}) & \rTo & 0.
\end{diagram}
\end{equation}
As it is well known \cite{KasparovJSM1981,JensenBook}, if $C_r^\ast(\Gg|_{\Xi_\infty^c})$ is stable, then such exact sequences are classified by a semi-group, which in many instances, becomes a group. This group then is isomorphic to Kasparov's group 
$$KK^1\big(C_r^\ast(\Gg|_{\Xi_\infty}),C_r^\ast(\Gg|_{\Xi_\infty^c})\big)$$ 
and the bulk-defect correspondence principle is contained in Kasparov product
\begin{equation}\label{Eq:KK0}
KK^\ast\big (\CM,C_r^\ast(\Gg|_{\Xi_\infty})\big ) \times KK^1\big(C_r^\ast(\Gg|_{\Xi_\infty}),C_r^\ast(\Gg|_{\Xi_\infty^c})\big) \to KK^\ast\big (\CM,C_r^\ast(\Gg|_{\Xi_\infty^c})\big ).
\end{equation}
Indeed, if the left side can be shown to be nontrivial for a projection/unitary associated with a Hamiltonian from the bulk algebra  $C_r^\ast(\Gg|_{\Xi_\infty})$, then a spectral statement can be made about arbitrary lifts of the bulk Hamiltonian into the full algebra $C_r^\ast(\Gg)$ (see \cite{BourneMPL2015} and section~\ref{Sec:SpSt}).

Thus, our work generalizes the $KK$-theoretic machinery developed for samples with boundaries \cite{BourneMPL2015,BourneAHP2017} to generic lattice defects. In fact, the new approach gives an exhaustive characterization of the lattice defects, in the following sense:

\begin{theorem} Any lattice defect can be canonically associated with an element of the Kasparov group $KK^1\big(C_r^\ast(\Gg|_{\Xi_\infty}),C_r^\ast(\Gg|_{\Xi_\infty^c})\big)$, via the mechanism described above. Furthermore, two lattice defects connected to the same Kasparov element have identical bulk-defect correspondences.
\end{theorem} 

The big picture that emerges, at least in the examples worked out here, is that of continuous fields of $C^\ast$-algebras, separately indexed by $\Xi_\infty^i$'s, which are glued together by the groupoid $C^\ast$-algebra $C_r^\ast(\Gg)$. The various way of gluing are classified by specific Kasparov groups, which depend on how we cluster the asymptotic manifolds $\Xi_\infty^i$ when deciding how to split the space of units. Furthermore, the defects can be nucleated and fused, and the rules for these processes can be described by the Kasparov groups we just mentioned.

To demonstrate the above ideas, we first present in Sec.~\ref{Sec:OldExample} a classical well understood example of a bulk-boundary correspondence principle, which we re-analyze here within the newly proposed framework. In Sec.~\ref{Sec:NewExample}, we apply the formalism to a 2-dimensional lattice with a disclination defect, for which we present numerical simulations of non-trivial bulk-defect correspondence principles. As we shall see, for a lattice model which supply the top generator of the first Kasparov group in Eq.~\ref{Eq:KK0}, there are topological states localized around a single point, the core of the defect, whose existence and characteristics are fully determined by the type of the defect and by the asymptotic properties of the lattice models. In materials science, the lattice defect mentioned above is known as a topological defect. They are rationalized by a topological charge given by the so called deficit angle (see Sec.~\ref{Sec:Discussion}). Using our formalism, we identify the unit of this topological charge with the top generator of the Kasparov group $KK^1 (C^\ast(\ZM^3),\CM)$.

Lastly, we mention that there is partial overlap between our numerical examples and some from existing literature, such as \cite{TeoPRL2013,BenalcazarPRB2014,VarjasPRL2019,GeierArxiv2020}. These works explore the role of crystalline/point group symmetries, while such symmetries do not play any role in our analysis. In particular, the bulk-defect correspondence principle established here is robust against any lattice deformation that does not change the topology of the space of units.

\section{The groupoid Algebra Associated to Linear Dynamics}
\label{Sec:DiscreteRes}

\subsection{Galilean invariant dynamical matrices}
\label{SubSec:LocalAlg0}

The setting is that of a network $\Ll$ of coupled {\it identical} quantum or classical resonators. We assume that the internal structure of the resonators has been fixed and that they couple to each other through a force field determined entirely by their internal structure. As a result, the coupling constants between the resonators are entirely determined by the pattern of their spatial arrangement. Very importantly, we assume that no background fields are present. In the linear dynamical regime, the dynamics of the coupled degrees of freedom is determined by a self-adjoint bounded operator over $\ell^2(\Ll)$, usually called the dynamical matrix,
\begin{equation}\label{Eq:DynMat}
D_\Ll = \sum_{x,x' \in \Ll} w_{x,x'}(\Ll) \, |x\rangle \langle  x' |, \quad w_{x,x'}(\Ll) \in \CM,
\end{equation}
with the sum converging in the strong operator norm. If the resonators have more than one internal degree of freedom, the coefficients $w_{x,x'}(\Ll)$ become matrices. Since this is almost a trivial difference, we proceed for now by assuming just one degree of freedom per resonator.

As explained above, once the internal structure of the resonators is fixed, there is a well defined correspondence $\Ll \mapsto D_\Ll$ between patterns and dynamical matrices, which in principle can be mapped in the laboratory. If two such patterns enter the relation $\Ll' = \Ll- y$ for some $y\in \RM^d$, then necessarily
\begin{equation}
w_{x,x'}(\Ll) = w_{x-y,x'-y}(\Ll - y), \quad x,x' \in \Ll,
\end{equation}
which is a consequence of the Galilean invariance of the governing physical laws. A direct consequence of these relations is a reduction of Eq.~\eqref{Eq:DynMat} to
\begin{equation}
D_\Ll = \sum_{x\in \Ll} \sum_{y\in \Ll-x} w_{0,y}(\Ll-x) \, |x \rangle \langle x+y |,
\end{equation}
or
\begin{equation}\label{Eq:CanonicalD}
(D_\Ll \psi)(x) = \sum_{y\in \Ll-x} w_{y}(\Ll-x) \, \psi(x+y), \quad \psi \in \ell^2(\Ll),
\end{equation}
where we dropped the trivial index. The expression in Eq.~\eqref{Eq:CanonicalD} is directly related to the regular representations of the groupoid algebra canonically attached to a pattern  \cite{Bellissard1986,Bellissard1995,KellendonkRMP95}, as explained next.

\subsection{The canonical groupoid associated with the dynamics}
\label{SubSec:Patterns}

When dealing with  patterns in $\RM^d$, a very useful metric space is that of compact subsets equipped with the Hausdorff metric, denoted here by $\big(\Kk(\RM^d),d_{\rm H}\big)$. When dealing with non-compact sub-sets, one needs to pass to the larger space $\Cc(\RM^d)$ of closed sub-sets of $\RM^d$. Since the latter can be canonically embedded in $\SM^d$, the $d$-dimensional sphere, one can use the Hausdorff metric on $\SM^d$ to equip $\Cc(\RM^d)$ with a metric, which we denote by $\bar d_{\rm H}$. This and other equivalent characterizations of $\big ( \Cc(\RM^d), \bar d_{\rm H}\big )$ can be found in \cite{ForrestAMS2002,LenzTheta2003}. We will refer to the metric space $\big ( \Cc(\RM^d),\bar d_{\rm H} \big )$ as the space of patterns in $\RM^d$. It is bounded, compact and complete. Furthermore, there is a continuous action of $\RM^d$ by translations, that is, a homomorphism between topological groups,
\begin{equation}\label{Eq:RDGroupAction}
\mathfrak t : \RM^d \rightarrow {\rm Homeo}\big ( \Cc(\RM^d),\bar d_{\rm H} \big ), \quad \mathfrak t_{x}(\Lambda) = \Lambda - x.
\end{equation} 

We will not work with the entire space of closed sets. Instead, we will focus on discrete sub-sets and, more specifically, on uniformly discrete sub-sets. We recall that $\Ll \in \Cc(\RM^d)$ is said to be $r$-uniformly discrete if any closed ball of radius $r \in \RM_+$ in $\RM^d$ contains at most one point \cite{Bellissard2000}. Let us point out that, the theory of homogeneous condensed matter systems is developed inside the class of Delone sets \cite{Bellissard2003}. However, when dealing with defects in materials, we are forced to step out that comfort zone and work with just uniformly discrete sets. The latter is still a manageable class for the space of $r$-uniformly  discrete patterns is a compact subset of $\Cc(\RM^d)$, invariant against translations \cite{Bellissard2000}. Furthermore, the notion of continuous hull and of its canonical transversal is still useful:

\begin{definition}[\cite{Bellissard2000}]\label{Def:Hull} The continuous hull of a uniformly discrete pattern $\Ll_0$ is the topological dynamical system $(\Omega_{\Ll_0},\mathfrak t,\RM^d)$, where
\begin{equation}
\Omega_{\Ll_0}= \overline{\{\mathfrak t_x(\Ll_0)=\Ll_0 - x, \ x\in \RM^d\}} \subset \Cc(\RM^d).
\end{equation}
The canonical transversal of $(\Omega_{\Ll_0},\mathfrak t,\RM^d)$  is defined as
\begin{equation}
\Xi_{\Ll_0} = \{ \Ll \in \Omega_{\Ll_0}, \ 0 \in \Ll\}.
\end{equation}
Both $\Omega_{\Ll_0}$ and $\Xi_{\Ll_0}$ are compact subspaces of $\Cc(\RM^d)$.
\end{definition}

Every point $\Ll \in \Omega_{\Ll_0}$ defines a closed subset of $\RM^d$. Theorem~2.8 of \cite{Bellissard2000} assures us that all these closed subsets are in fact uniformly discrete subsets. Note, however, that, if $\Ll \in \Omega_{\Ll_0}$ and $\Ll$ does not belong to the orbit of $\Ll_0$, then the orbit of $\Ll$ might not be dense in $\Omega_{\Ll_0}$. An example which displays such phenomenon is a periodic pattern in $\RM$ with one defect \cite{SadunBook}[p.~9]. Another example is the disordered pattern $\{n + \xi_n\}_{n \in \ZM}$, with $\xi_n$ drawn randomly from some interval. In both cases, the periodic pattern $\ZM$ sits inside the hull and the orbit of this pattern traces a circle, which is strictly smaller than the hull of the original pattern. In fact, the character of the patterns contained inside the hull can vary drastically. Indeed, out of two un-related Delone sets from $\RM$, one can generate a third one by joining the left side of one pattern with the right side of the other one. The resulting pattern will have a hull that contains the hulls of the original two patterns.\footnote{Assuming that one-sided translations generate the whole hull of the original patterns.} Of course, both notions introduced in \ref{Def:Hull} become useful in applications if they can be explicitly computed. 

We now introduce the Bellissard-Kellendonk groupoid associated to a uniformly discrete pattern \cite{Bellissard1986,Bellissard1995,KellendonkRMP95}: 

\begin{definition} The topological groupoid of an $r$-uniformly discrete pattern $\Ll_0$ consists of:
\begin{enumerate}[\rm 1.]

\item The set
\begin{equation}
\Gg :=\{(\Ll,x) \in \Xi_{\Ll_0} \times \RM^d, \ x \in \Ll\}
\end{equation}
equipped with the inversion map
\begin{equation}
(\Ll,x)^{-1} = (\Ll-x,-x).
\end{equation}
\item The sub-set of $\Gg \times \Gg$
\begin{equation}
 \Gg^{(2)} =\{((\Ll,x),(\Ll',y)) \in \Gg \times \Gg, \  \Ll'=\Ll-x\}
 \end{equation}
equipped with the composition 
\begin{equation}
(\Ll,x) \cdot (\Ll-x,y) = (\Ll,x+y).
\end{equation}
\end{enumerate} 
The topology on $\Gg$ is the relative topology inherited from $\Xi_{\Ll_0} \times \RM^d$.
\end{definition}

\begin{proposition}[\cite{BourneAHP2020},~Prop.~2.16] $\Gg$ is a second-countable, Hausdorff, \'etale groupoid.
\end{proposition}

In the world of groupoids, the \'etale groupoids are the equivalent of discrete groups from the world of locally compact groups. Indeed, the fibers are always discrete topological spaces and the Haar system is supplied by the counting measures \cite{SimsSzaboWilliamsBook2020}[Ch.~8]. In our specific case, the source and range maps are
\begin{equation}
\begin{aligned}
& \mathfrak s(\Ll,x) := (\Ll,x)\cdot (\Ll,x)^{-1} = (\Ll-x,0), \\ 
& \mathfrak r(\Ll,x) := (\Ll,x)^{-1}(\Ll,x)= (\Ll,0),
\end{aligned}
\end{equation}
hence
\begin{equation}\label{Eq:SmRm}
\mathfrak r^{-1}(\Ll)=\{(\Ll,x), \, x \in \Ll\}, \quad \mathfrak s^{-1}(\Ll)=\{(\Ll-x,-x), \, x \in \Ll\},
\end{equation}
and $\mathfrak r^{-1}(\Ll)$ can be canonically identified with $\Ll$. Also, since the inversion map takes $\mathfrak r^{-1}(\Ll)$ into $\mathfrak s^{-1}(\Ll)$, the space $\mathfrak s^{-1}(\Ll)$ can be also canonically identified with $\Ll$. We recall that all $\Ll \in \Xi$ are $r$-uniformly discrete point sets. Furthermore, the space of units $\Gg^{(0)}$, {\it i.e.} the image of $\Gg$ through $\mathfrak s$ or $\mathfrak r$, can be canonically identified with $\Xi_{\Ll_0}$.

\subsection{The canonical $C^\ast$-algebra associated with the dynamics}

For each \'etale groupoid, there are two canonically associated $C^\ast$-algebras \cite{RenaultBook}. We are interested in the reduced $C^\ast$-algebra $C_r^\ast(\Gg)$, which contains the space of compactly supported functions $f: \Gg \to \CM$  as a dense sub-algebra and the product of such elements is
\begin{equation}\label{Eq:Multi}
(f \cdot f')(\Ll,x) = \sum_{y \in \Ll} f(\Ll,y) f'(\Ll-y,x-y).
\end{equation}
Also, the involution takes the form $f^\ast(\Ll,x) = \bar f (\Ll -x,-x)$. Its left regular representations are indexed by the unit space, hence by $\Ll \in \Xi_{\Ll_0}$, and take place on the Hilbert spaces $\ell^2(s^{-1}(\Ll,0))$. We use the map  $\mathfrak j(\Ll-x,-x)=x$ to establish a canonical isomporphism between the Hilbert spaces $\ell^2(s^{-1}(\Ll,0))$ and $\ell^2(\Ll)$ and write the representations as
\begin{equation}
[\pi_{\Ll}(f)\psi](x) = (f\cdot (\psi\circ \mathfrak j))(\mathfrak j^{-1}(x)), \quad \psi \in \ell^2(\Ll), \quad x \in \Ll,
\end{equation}
or, in a more explicit and convenient form,
\begin{equation}\label{Eq:Rep1}
[\pi_{\Ll}(f)\psi](x) = \sum_{y \in \Ll-x} f(\Ll-x,y) \psi(x+y).
\end{equation}
If we compare the expression \eqref{Eq:Rep1} with the action of Galilean invariant dynamical matrices \eqref{Eq:CanonicalD}, we see that they are identical once we identify $f(\Ll-x,y)$ and $w_y(\Ll-x)$. The conclusion is that all Galilean invariant dynamical matrices can be generated using the regular representations of $C_r^\ast(\Gg)$, which is the smallest algebra with this property. If the resonators have more than one internal degree of freedom, then one passes to the stabilization $\KM(\NM) \times C_r^\ast(\Gg)$, where $\KM(\NM)$ is the algebra of compact operators over $\ell^2(\NM)$.

\section{Bulk-boundary correspondence via groupoid methods}

The bulk-boundary-defect correspondence principle is at the core of any theory of topological phases of matter. In \cite{KellendonkRMP2002},  Kellendonk {\it et al} identified the engine of the bulk-boundary principle to be a certain exact sequence of $C^\ast$-algebras. We review here these general principles, this time in the context of groupoid $C^\ast$-algebras.

\subsection{Invariant ideals}

We first comment on the notion of $\Gg$-invariant sub-sets in our particular context. Let $\Sigma$ be a subset of the space of units $\Xi$ of $\Gg$. Then $\Sigma$ is said to be $\Gg$-invariant if $\mathfrak s^{-1}(\Sigma) = \mathfrak r^{-1}(\Sigma)$. In our particular context, this means that, if $\Ll \in \Sigma$, then $\Ll -x \in \Sigma$ for all $x \in \Ll$. Suppose now that $\Sigma$ is an open $\Gg$-invariant subset. Then \cite{WilliamsBook}[Sec.~5.1]
\begin{equation}
\Gg|_\Sigma := \mathfrak s^{-1}(\Sigma) = \mathfrak r^{-1}(\Sigma)
\end{equation} 
is an \'etale groupoid. Furthermore, $C^\ast(\Gg|_\Sigma)$
is a two-sided ideal in $C^\ast_r(\Gg)$. On the other hand, the closed complement $\Sigma^c = \Xi \setminus \Sigma$ is also a $\Gg$-invariant sub-set and there exists the short exact sequence of $C^\ast$-algebras \cite{WilliamsBook}[Th.~5.1]
\begin{equation}\label{Eq:ExactSeq1}
\begin{diagram}
0 & \rTo & C_r^\ast(\Gg|_\Sigma) & \rTo{\rm i} & C_r^\ast(\Gg) & \rTo{\rm j} & C_r^\ast(\Gg|_{\Sigma^c}) & \rTo & 0.
\end{diagram}
\end{equation}
It is important to remember that such an exact sequence always exists for the full groupoid $C^\ast$-algebras and that it always exists for the reduced groupoid $C^\ast$-algebra if the reduced and universal groupoid $C^\ast$-norms coincide \cite{WilliamsBook}[Prop.~5.2]. If the reduced and universal norms do not coincide, there are explicit examples where the sequence of reduced groupoid algebras is not exact \cite{HigsonGAFA2002}. Surprisingly, the counterexamples found in \cite{HigsonGAFA2002} are related in spirit to our applications, with the crucial difference that the groups appearing in our examples are amenable. Therefore, one should take extra care if our framework is to be applied to hyperbolic or fractal crystals, where the underlying groups are not amenable.

\subsection{Spectral statements via $K$-theory} 
\label{Sec:SpSt}

Any short exact sequence of $C^\ast$-algebra induces a canonical 6-term exact sequence for the $K$-theories of the algebras. In particular, for the exact sequence~\eqref{Eq:ExactSeq1}, we have
\begin{equation}\label{Eq-SixTermDiagram}
\begin{diagram}
& K_0\big(C_r^\ast(\Gg|_\Sigma)\big ) & \rTo{ \ {\rm i}_\ast \ \ } & K_0\big(C_r^\ast(\Gg)\big) & \rTo{\ \ {\rm j}_\ast \ \ } & K_0\big(C_r^\ast(\Gg|_{\Sigma^c})\big) &\\
& \uTo{\rm Ind} & \  &  \ & \ & \dTo{\rm Exp} & \\
& K_1\big(C_r^\ast(\Gg|_{\Sigma^c})\big)  & \lTo{ \ \ {\rm j}_\ast} & K_1\big(C_r^\ast(\Gg)\big) & \lTo{\ \ {\rm i}_\ast \ } & K_1\big(C_r^\ast(\Gg|_{\Sigma})\big) &
\end{diagram}
\end{equation}

The general principle behind the spectral statements derived from the above 6-term exact sequence is as follows. Suppose $h$ is a self-adjoint element from $C_r^\ast(\Gg|_{\Sigma^c})$. A gap in the spectrum ${\rm Spec}(h)$ of $h$ is a connected component of $\RM \setminus {\rm Spec}(h)$. Let $G$ be such a gap and let $\hat h$ be any lift of $h$ in $C_r^\ast(\Gg)$. Note that, necessarily, ${\rm Spec}(h) \subseteq {\rm Spec}(\hat h)$. Then, if $p_G$ is the spectral projection of $h$ onto the spectrum below $G$ and if ${\rm Exp}[p_G]_0 \neq [1]_1$, then necessarily $G \subset {\rm Spec}(\hat h)$ (see \cite{ProdanJGP2019} for details).

We can also use the left side of the diagram to generate a spectral statement. Indeed, let $h$ be self-adjoint element from $M_{2N}(\CM) \otimes C_r^\ast(\Gg|_{\Sigma^c})$, where $M_k(\CM)$ denotes the algebra of $k \times k$ matrices with complex entries. We denote the unit of this algebra by $1_k$. We assume that $h$ has the chiral symmetry, more specifically, that
\begin{equation}
J \, h \, J^{-1} = -h, \quad J = \begin{pmatrix} 1_N & 0 \\ 0 & -1_N \end{pmatrix}\otimes {\rm id}.
\end{equation}
This property automatically implies that ${\rm Spec}(h)$ is symmetric relative to the origin. Next, we assume that $h$ has a gap $G$ in its spectrum such that $0 \in G$. 
 Let $\hat h$ be any lift of $h$ in $M_{2N}(\CM) \otimes C_r^\ast(\Gg)$ such that $\hat h$ also display the chiral symmetry. Now, let $p_G$ be the spectral projection of $h$ onto the spectrum below $G$. Due to the chiral symmetry, $p_G$ takes the following particular form
\begin{equation}
p_G = \frac{1}{2} \begin{pmatrix} 1_N & - u_G^\ast \\ - u_G & 1_N \end{pmatrix}
\end{equation}
with $u_G$ a unitary element from $M_{N}(\CM) \otimes C_r^\ast(\Gg|_{\Sigma^c})$. Then, if ${\rm Ind}[u_G]_1 \neq [0]_0$, then necessarily $0 \in {\rm Spec}(\hat h)$ \cite{ProdanSpringer2016}[Corollary~4.3.4]. Since this is  used in our examples in Sec.~\ref{Sec:NewExample}, let us be more concrete and specify a convenient presentation of the index map \cite{ProdanSpringer2016}[Sec.~4.3.2]:
\begin{equation}
{\rm Ind}[u_G]_1 = [e_G]_0 - [{\rm diag}(1_N,0_N)]_0
\end{equation}
where
\begin{equation}
\begin{aligned}
e_G & = e^{-\imath \frac{\pi}{2} g(\hat h)} {\rm diag}(1_N,0_N) e^{\imath \frac{\pi}{2} g(\hat h)} \\
& = \tfrac{1}{2} J(e^{\imath \pi g(\hat h)}-1_{2N}) + {\rm diag}(0_N,1_N),
\end{aligned}
\end{equation} 
with $g: \RM \to \RM$ {\it any} odd function equal to $\pm 1$ above/below the $G$. From the above expression, one can see explicitly that, if $\hat h$ is gapped at 0, then the variation of $g$ can be concentrated entirely inside that gap and, as a consequence, the index map returns a trivial value. In special cases, $e_G$ can be connected to a spectral projection of $\hat h$. Indeed, if $\hat h$ displays spectral gaps above and below 0 and these spectral gaps are located inside $G$, then the variation of $g$ can be concentrated inside these gaps and, as a consequence, $e_G = J\chi_{[-\delta,\delta]}(\hat h)+{\rm diag}(0_N,1_N)$, where $\pm \delta$ are the midpoints of the spectral gaps of $\hat h$ we just mentioned. This, for example, applies to the model presented in Sec.~\ref{Sec:OldExample} but it does not apply to the model presented in Sec.~\ref{Sec:NewExample}.

It is useful to place the above discussion in an even more general context, namely, that of Kasparov's $KK$-theory \cite{KasparovJSM1987,JensenBook}, following \cite{BourneMPL2015} as a model. Specifically, the extension~\eqref{Eq:ExactSeq1} defines an element $[{\rm ext}]^1$ of the group $KK^1\big (C_r^\ast(\Gg|_{\Sigma^c}),C_r^\ast(\Gg|_{\Sigma}\big )$ and the connecting maps in the 6-term diagram~\eqref{Eq-SixTermDiagram} can be expressed as Kasparov products \cite{BlackadarBook}[Th.~19.5.7]:
\begin{equation}
{\rm Exp}[p_G]_0 = [p_G]_0 \times [{\rm ext}]^1 \in KK^1\big (\CM,C_r^\ast(\Gg|_{\Sigma})\big ) \simeq K_1 \big (C_r^\ast(\Gg|_{\Sigma}) \big ),
\end{equation}
and 
\begin{equation}
{\rm Ind}[u_G]_1 = [u_G]_1 \times [{\rm ext}]^1 \in KK^0\big (\CM,C_r^\ast(\Gg|_{\Sigma})\big ) \simeq K_0 \big (C_r^\ast(\Gg|_{\Sigma}) \big ).
\end{equation}
As we shall see, different lattice defects correspond to different elements of the extension group. For practical applications, one needs to compute the $K$-theories of both $C_r^\ast(\Gg|_{\Sigma^c})$ and $C_r^\ast(\Gg|_{\Sigma})$ algebras, as well as the action of the connecting maps on the generators of the $K$-groups. Such computations have been carried out, for example, in \cite{ProdanSpringer2016} for crossed product $C^\ast$-algebras by $\ZM^d$ and, as a result, the spectral statements in \cite{ProdanSpringer2016} are formulated in more precise terms then it was done above. Same calculations combined with the analysis from \cite{ProdanRMP2016} will prove instrumental for the example discussed in Sec.~\ref{Sec:NewExample}.

Lastly, the reader should be aware that the bulk-boundary correspondence contains an additional dynamical statement, which refers to the localized or de-localized character of the spectrum of $\hat h$ inside $G$ [see Sec.~6.6 in \cite{ProdanSpringer2016}]. The dynamical statements were established in \cite{ProdanSpringer2016} via index theory and we will not touch this subject in the present work. 

\section{Examples of bulk-boundary-defect correspondences}

\subsection{An old example analyzed in the new framework}
\label{Sec:OldExample} 

We consider here the minimal model that displays a non-trivial bulk-boundary correspondence. It has been thoroughly analyzed in Sec.~1 of \cite{ProdanSpringer2016}. As it is well known, the short exact sequence of $C^\ast$-algebras behind its bulk-boundary correspondence mechanism is isomorphic to the classical Toeplitz extension. Our goal here is to show how same extension can be formulated with groupoid methods, specifically, as in Eq.~\eqref{Eq:ExactSeq1}.

On the infinite lattice, the model is defined on the Hilbert space $\CM^{2N} \otimes \ell^2 (\ZM)$ and takes the form
\begin{equation}\label{Eq:BulkModel1}
H = \sum_{q \in \ZM} W_q \otimes S^q, \quad W_q = W_{-q}^\ast \in M_{2N}(\CM),
\end{equation}
where $S$ is the shift operator on $\ell^2(\ZM)$, $S|n\rangle = |n-1\rangle$. Hence, the model comes from the stabilization of the algebra $C^\ast(S)$ generated by a single unitary element. When the models are restricted to half-space, {\it i.e.} to $\ell^2(\NM)$, it is natural to consider the stabilization of the algebra $C^\ast(\hat S)$, generated by the half-shift on $\ell^2(\NM)$. We recall that the half-shift and its adjoint enter the following relations
\begin{equation}
\hat S \, \hat S^\ast = 1, \quad \hat S^\ast \, \hat S = 1- E, \quad E = |0\rangle \langle 0 |.
\end{equation}
Let $C^\ast(\hat S) \, E \, C^\ast(\hat S)$ be the closed two-sided ideal generated by the projection $E$. It is not difficult to see that this ideal contains operators with matrix elements decaying to zero away from the boundary. Then we have the exact sequence of $C^\ast$-algebras
\begin{equation}\label{Eq:Seq1}
\begin{diagram}
0 & \rTo &C^\ast(\hat S)\, E \, C^\ast(\hat S) &\rTo{\ \ {\rm i} \ \ }  & C^\ast(\hat S)  &\pile{ \rTo{\mathrm{j}} \\ \lTo_{\mathrm{i}'} } & C^\ast(S) &\rTo &0,
\end{diagram}
\end{equation}
which is split at the level of linear spaces and isomorphic to the Toeplitz extension [see Sec.~3.2.3 in \cite{ProdanSpringer2016}].

Now, let $H$ be as in Eq.~\eqref{Eq:BulkModel1}. Any self-adjoint lift $\hat H \in M_{2N}(\CM) \otimes  C^\ast(\hat S)$ of $H$ generates a model with a boundary for $H$. Indeed, the linear splitting of the exact sequence~\eqref{Eq:Seq1} assures us that any such lift can be uniquely written as $
\hat H = V^\ast H V + \tilde H$, where $V$ is the isometry from $\ell^2(\NM)$ to $\ell^2(\ZM)$ and $\tilde H$ is from $C^\ast(\hat S) \, E \, C^\ast(\hat S)$. This decomposition says that $\hat H$ is generated from $H$ by imposing Dirichlet boundary condition, followed by addition of a potential $\tilde H$ localized near the boundary. Since this recipe generates all physically sound boundary models for $H$, the stabilization of the algebra $C^\ast(\hat S)$ can be rightfully called the algebra of models with a boundary. Furthermore, the stabilization of the algebra $C^\ast(\hat S)\, E \, C^\ast(\hat S)$ can be rightfully called the algebra of boundary potentials. 

Due to its simplicity, much is known about this example. Specifically, the algebra of boundary potentials is isomorphic to $\KM\big(\ell^2(\NM)\big )$ and $[E]_0$ can be chosen as the generator of its $K_0$-group. The $K_1$-group of $C^\ast(S)$ is generated by $[S]_1$ and the action of the index map is ${\rm Ind}[S]_1 =[E]_0$. The bulk-boundary principle is also very straightforward: If $H$ displays the chiral symmetry and the class of $u_G$ in the $K_1$-group accepts a decomposition $[u_G]_1 = n[S]_1$ with $n\neq 0$, then $0 \in {\rm Spec}(\hat H)$ for any lift $\hat H$ with chiral symmetry. In fact, one can say a bit more, specifically, that the degeneracy of the spectrum at $0$ is a least $n$. A simple example where these statements can be verified by explicit calculations is supplied by the Hamiltonians 
\begin{equation}
H = \begin{pmatrix} 0 & S^\ast \\ S & 0 \end{pmatrix}, \quad \hat H = \begin{pmatrix} 0 & \hat S^\ast \\ \hat S & 0 \end{pmatrix},
\end{equation} 
for which one will find ${\rm Spec}(H) = \{-1,1\}$ and ${\rm Spec}(\hat H) = \{-1,0,1\}$. The above analysis assures us that the zero eigenvalue cannot be removed by any boundary potential with chiral symmetry.

\begin{figure}[t]
\center
\includegraphics[width=0.8\textwidth]{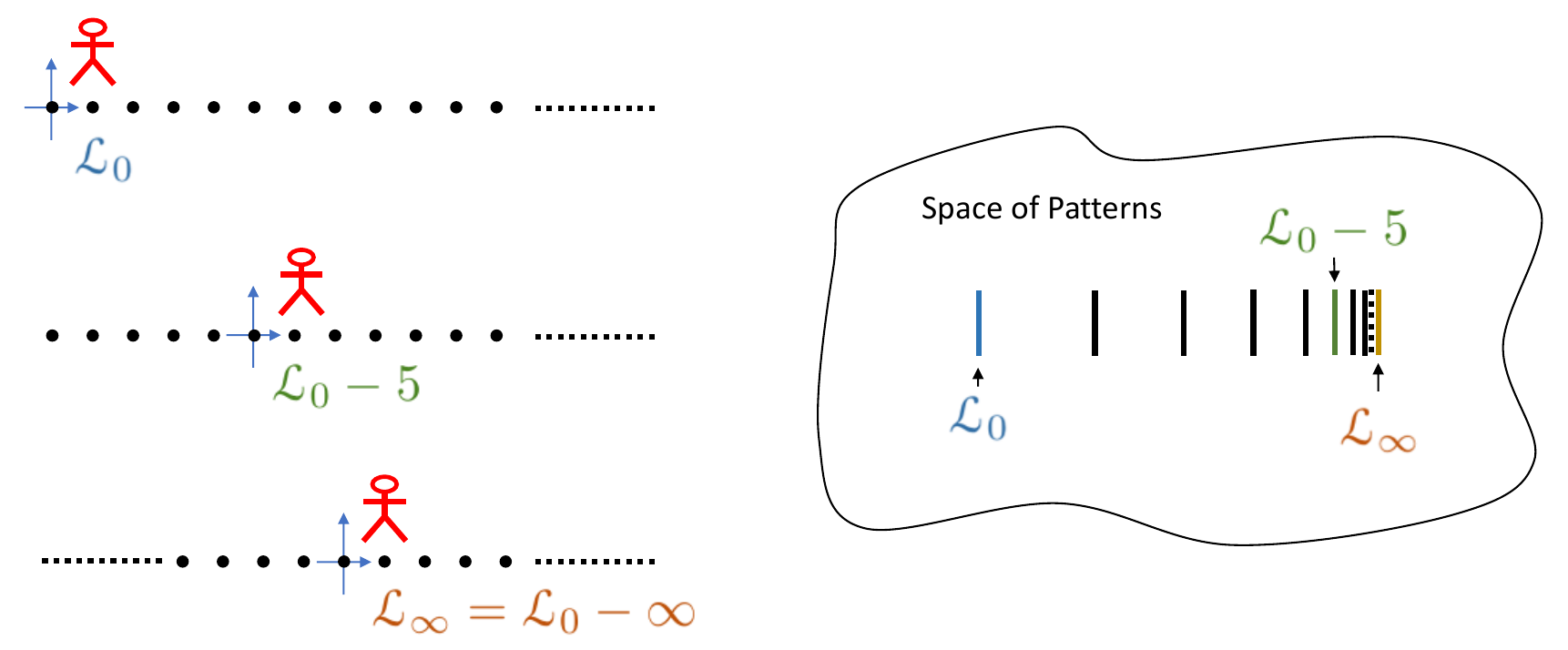}\\
  \caption{\small (Left) An observer placed at the origin assesses translated half-lattices. (Righ) Qualitative representation of the sequence of translated half-lattices $\Ll_0-n$ in the space of patterns $\Cc(\RM)$. Just for visualization, the half-lattices in the space of patterns are shown as segments rather than points.
  }
 \label{Fig:Transv1}
\end{figure}

We now explain how all those phenomena can be seen within the formalism based on the groupoid algebra. As opposed to \cite{BourneAHP2020}, our analysis starts from the half-space models. Henceforth, we consider the half-lattice represented by $\Ll_0=\NM$. Our first task is to compute the transversal of this uniformly discrete pattern and we will use Fig.~\ref{Fig:Transv1} for this. As a guiding principle, the distance $\bar d_{\rm H}$ between two patterns becomes smaller and smaller as the patterns resemble each other on larger and larger neighborhoods of the origin. To the observer sitting at $x=0$, the patterns $\Ll_0$, $\Ll_0-1$ will appear quite different from each other, hence the distance between them is large, in fact, equal to one. The lattices $\Ll_0-5$ and $\Ll_0-6$ will appear the same over the neighborhood $[-5,5]$, hence the distance between them will be much smaller than 1. One can then see that the distance $\bar d_{\rm H}$ between consecutive patterns $\Ll_0-n$ and $\Ll_0-n-1$ decreases with $n$ and, as such,  the sequence of lattices $\Ll_0-n$ and the transversal $\Xi_{\Ll_0}$ must look as in Fig.~\ref{Fig:Transv1}. The latter is a sequence of discrete points with one accumulation point. This accumulation point is included because, by definition, the transversal is the closure of the discrete orbit of $\Ll_0$. It corresponds to the pattern $\Ll_\infty = \ZM$ and it has the distinct property that $\Ll_\infty - x = \Ll_{\infty}$ for all $x \in \Ll_\infty$. This property implies that [see Eq.~\eqref{Eq:SmRm}]
\begin{equation}
r^{-1}(\Ll_\infty) = s^{-1}(\Ll_\infty)=\ZM,
\end{equation} 
hence $\Ll_\infty$ is a closed $\Gg$-invariant set. Automatically, the open complement $\Sigma = \Xi_{\Ll_0} \setminus \Ll_\infty$ is a $\Gg$-invariant set and we have the exact sequence of groupoid $C^\ast$-algebras
\begin{equation}\label{Eq:ExactSeq2}
\begin{diagram}
0 & \rTo & C_r^\ast(\Gg|_\Sigma) & \rTo{\rm i} & C_r^\ast(\Gg) & \rTo{\rm j} & C_r^\ast(\Gg|_{\Ll_\infty}) & \rTo & 0.
\end{diagram}
\end{equation}

We now establish an explicit isomorphism between the exact sequences~\eqref{Eq:ExactSeq1} and \eqref{Eq:ExactSeq2}. For this we define the element $s \in C_r^\ast(\Gg)$,
\begin{equation}
s(\Ll,x) = \chi_{\Ll}(1) \delta_{x,1},
\end{equation}
with the conjugate
\begin{equation}
s^\ast(\Ll,x) = s(\Ll-x,-x) = \chi_{\Ll}(-1) \delta_{x,-1}.
\end{equation}
Here, $\chi_\Ll$ is the indicator function of the set $\Ll \subset \RM$.

\begin{proposition} The two elements enter the following relations:
\begin{equation}
s \cdot s^\ast = 1, \quad s^\ast \cdot s = 1- e, \quad e=e^2 = e^\ast,
\end{equation}
with $e$ given explicitly by $e(\Ll,x) = \delta_{\Ll,\Ll_0} \delta_{x,0}$.
\end{proposition}

\proof Using the multiplication rule~\eqref{Eq:Multi}, one can see that
\begin{equation}
s \cdot s^\ast(\Ll,x) = \delta_{x,0}, \quad s^\ast \cdot s (\Ll,x) = (1- \delta_{\Ll,\Ll_0}) \delta_{x,0}.
\end{equation}
The last term is exactly $1-e$.\qed

\vspace{0.2cm}

Our goal now is to show that $C^\ast(\Gg) = C^\ast(s, s^\ast)$ and we will accomplished this in two steps.

\begin{proposition}\label{Pro:X1} For any $f,g \in C_c(\Gg)$,
\begin{equation}\label{Eq:Id1}
f \cdot e \cdot g = \sum_{m,n\in \NM} f(\Ll_m,-m)g(\Ll_0,n) \ (s^\ast)^m \cdot e \cdot s^m.
\end{equation}
\end{proposition}

\proof Firstly, let us point out that $s \cdot e = e \cdot s^\ast =0$, hence monomial of the form $s^m \cdot e \cdot (s^\ast)^m$ with $m>0$ or $n>0$ are all identically zero. Secondly, let us introduce the notation $\Ll_m = \Ll_0 - m$, $m\in \NM$, and $x_\Ll^0$ for the most left point of $\Ll \in \Sigma$. Now, using the rule for multiplication, one finds
\begin{equation}
(f \cdot e \cdot g)(\Ll,x) = f(\Ll,x_\Ll^0) g(\Ll_0,x-x_\Ll^0) \chi_\Sigma(\Ll),
\end{equation}
while
\begin{equation}
[(s^\ast)^m \cdot e \cdot s^n](\Ll,x) = \delta_{\Ll,\Ll_m} \delta_{x,n-m}.
\end{equation}
Then Eq.~\eqref{Eq:Id1} follows. \qed

\begin{proposition} The exact sequences~\eqref{Eq:ExactSeq1} and \eqref{Eq:ExactSeq2} are isomorphic.
\end{proposition}

\proof The principal ideal $C^\ast(s,s^\ast) \, e \, C^\ast(s,s^\ast)$ is spanned by elements of the form $\sum_{m,n\in \NM} a_{m,n} \ (s^\ast)^m \cdot e \cdot s^m$. As such, Proposition~\ref{Pro:X1} established the equality between the principal ideals 
\begin{equation}
C^\ast(\Gg) \, e \, C^\ast(\Gg) = C^\ast(s,s^\ast) \, e \, C^\ast(s,s^\ast).
\end{equation}
Also, the following isomorphisms are obvious
\begin{equation}
C^\ast(\Gg|_{\Ll_\infty}) \simeq C^\ast(\ZM) \simeq C^\ast(s,s^\ast)/C^\ast(s,s^\ast) \, e \, C^\ast(s,s^\ast).
\end{equation}
These together with the fact that $C^\ast(s,s^\ast)$ is already a sub-algebra of $C^\ast(\Gg)$ assures us that in fact we have the equality $C^\ast(s,s^\ast)=C^\ast(\Gg)$. Then the mapping $\hat S \leftrightarrow s$ establishes the isomorphism between the exact sequences~\eqref{Eq:ExactSeq1} and \eqref{Eq:ExactSeq2}.\qed

\vspace{0.2cm} The conclusion is that the bulk-boundary correspondence can be established entirely within the groupoid framework via the exact sequence~\eqref{Eq:ExactSeq2} of groupoid $C^\ast$-algebras.

\subsection{Bulk-boundary correspondence for a lattice defect} 
\label{Sec:NewExample}

\begin{figure}[t]
\center
\includegraphics[width=\textwidth]{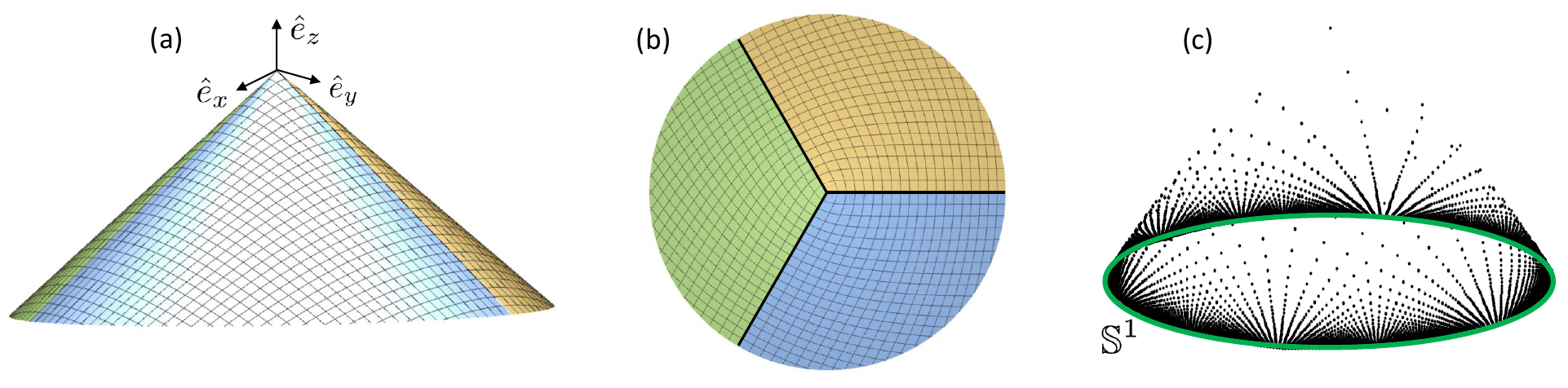}\\
  \caption{\small A lattice defect created by removing a quarter of a square lattice and gluing the exposed boundaries. Panels (a/b) show a view from the side/top of the resulting pattern $\Ll_0$, while panel (c) shows a representation of the transversal $\Xi_{\Ll_0}$. Panel (a) shows a coordinate system and its unit vectors used in the text.
  }
 \label{Fig:ConeNet}
\end{figure}

We analyze here the lattice defect described in Fig.~\ref{Fig:ConeNet}. The points of the resulted pattern $\Ll_0$ can be generated from $\ZM^2$ via the map
\begin{equation}\label{Eq:PattMap}
\ZM^2 \ni (n,m) \mapsto \Big (r \alpha \cos(\theta/\alpha), r \alpha \sin(\theta/\alpha),-r\sqrt{1-\alpha^2}\Big ) \in \RM^3,
\end{equation}
where $(r,\theta)$ are the polar coordinates of $(n,m)$ and $\alpha = 3/4$. For reasons that will become obvious below, we introduce the parameter $\xi$ such that $\alpha = \sin(\xi)$ and $\sqrt{1 -\alpha^2} = \cos(\xi)$. If we restrict $\theta$ to the interval $[0,\tfrac{3\pi}{2})$, then the map~\eqref{Eq:PattMap} is one to one. Obviously, any pair of points of $\ZM^2$ with same $r$ but $\theta=0$ and $\theta=\tfrac{3\pi}{2}$ are mapped into the same point of $\Ll_0$. Note that $(0,0)$ is mapped at the origin of $\RM^3$.

As before, the starting point for our analysis is the computation of the transversal. For this, we imagine the observer located at the origin of the physical space and looking at the translated pattern $\Ll_0 -x_0$ for some $x_0 \in \Ll_0$. Then, for every such $x_0$, the pattern will look distinct because the tip of the cone is in the field of view of the observer and the position of this tip can be observed to change with $x_0$. We now imagine shifting the pattern such that the tip of the cone is infinitely far away from the observer. In this case, the pattern will appear to the observer as a flat square lattice. To see this explicitly, let $(n_0,m_0) \in \ZM^2$ be the point that is mapped onto $x_0$ by~\eqref{Eq:PattMap}. Then a calculation based on \eqref{Eq:PattMap} will show that $\Ll_0 - x_0$ is given by the map
\begin{equation}
\ZM^2 \ni (n,m) \mapsto n \hat a_1+m \hat a_2 + o(r/r_0) \in \RM^3,
\end{equation}
where 
\begin{equation}
\begin{aligned}
\hat a_1(\theta_0)  & = \Big [\sin(\xi)\cos(\theta_0)\cos\big(\tfrac{\theta_0}{\alpha}\big)+\sin(\theta_0)\sin\big(\tfrac{\theta_0}{\alpha}\big)\Big] \hat e_x \\
&  \quad + \Big [\sin(\xi)\cos (\theta_0) \sin\big(\tfrac{\theta_0}{\alpha}\big) -\sin(\theta_0) \cos\big(\tfrac{\theta_0}{\alpha}\big)\Big ] \hat e_y \\
& \quad \quad - \cos(\xi)\cos(\theta_0) \hat e_z,
\end{aligned}
\end{equation}
and
\begin{equation}
\begin{aligned}
\hat a_2(\theta_0) & = \Big [\sin(\xi)\sin(\theta_0)\cos\big(\tfrac{\theta_0}{\alpha}\big)-\cos(\theta_0)\sin\big(\tfrac{\theta_0}{\alpha}\big)\Big] \hat e_x \\
& \quad + \Big [\sin(\xi)\sin (\theta_0) \sin\big(\tfrac{\theta_0}{\alpha}\big) + \cos(\theta_0) \cos\big(\tfrac{\theta_0}{\alpha}\big)\Big ] \hat e_y \\
& \quad \quad - \cos(\xi)\sin(\theta_0) \hat e_z.
\end{aligned}
\end{equation}
Above, $(r_0,\theta_0)$ are the polar coordinates of $(n_0,m_0)$ and $r = \sqrt{n^2 + m^2}$. It turns out that both $\hat a_1$ and $\hat a_2$ have norm one and they are orthogonal. As such, they can be generated from $\hat e_x$ and $\hat e_y$ via a 3-dimensional rotation $\widehat R(\theta_0)$ that is fully determined by $\theta_0$. This rotation can be easily read off from the expressions supplied above but its explicit form is not really needed. The conclusion is that, in the asymptotic limit $r_0 \rightarrow \infty$, the observer will see a square lattice $\Ll_\infty(\theta_0)$, which is just the lattice $\{n\hat e_x + m \hat e_y \in \RM^3, \ (n,m) \in \ZM^2\}$ rotated by $\widehat R(\theta_0)$. 

\begin{figure}[t]
\center
\includegraphics[width=\textwidth]{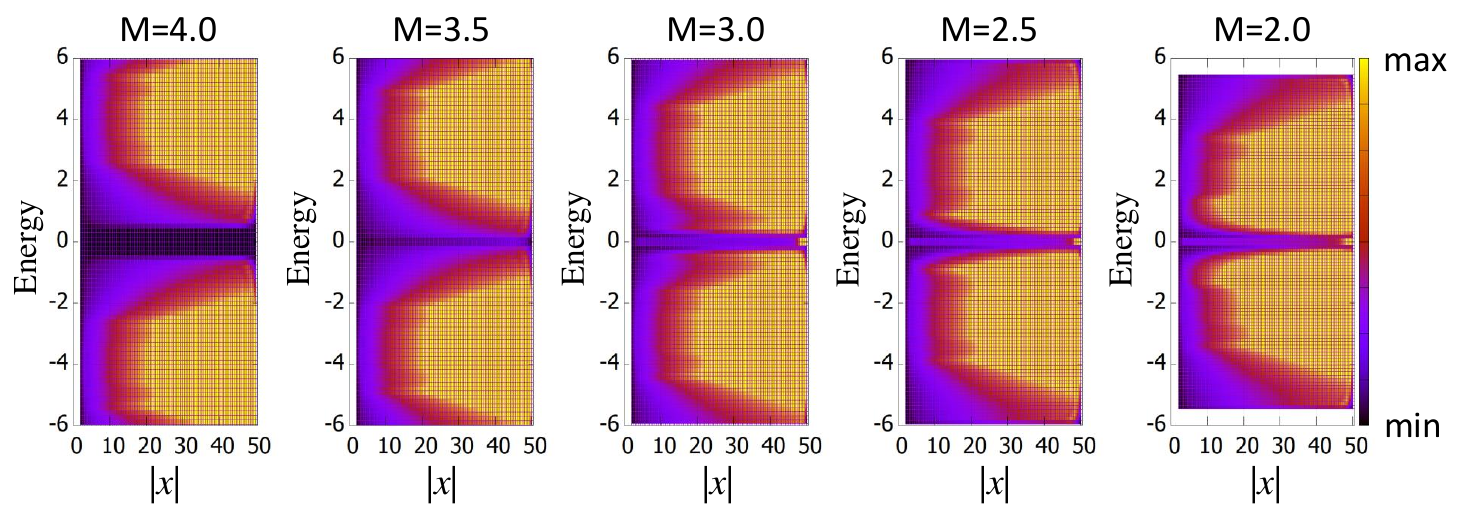}\\
  \caption{\small Evolution of the local density (LDOS) of states with the parameter $M$. LDOS is displayed as function of $R$, the distance to the tip of the cone. The notable feature is the presence/absence of zero-energy spectrum for various values of $M$, in perfect agreement with the topological phase transition discussed in the text. The calculation was performed with $R_{\rm max} = 50$.
  }
 \label{Fig:RadLdosVsM}
\end{figure}

An important observation is that
\begin{equation}\label{Eq:Twist}
\hat a_1(3\pi/2) = \hat a_2(0), \quad \hat a_2(3\pi/2) = -\hat a_1(0).
\end{equation}
This points to a nontrivial global twist of the square lattice experienced by the observer when walking around the rim of pattern shown in Fig.~\ref{Fig:ConeNet}(a). It also assures us that $\Ll_{\infty}(3\pi/2) = \Ll_{\infty}(0)$ and, as a result, the transversal $\Xi_{\Ll_0}$ looks as in Fig.~\ref{Fig:ConeNet}(c). As before, there is a natural decomposition $\Xi_{\Ll_0} = \Sigma \cup \Xi_\infty$, where $\Xi_\infty \simeq \SM^1$ and $\Sigma$ is the open set $\Xi_{\Ll_0} \setminus \Xi_\infty$. This prompts the 6-term exact sequence~\eqref{Eq-SixTermDiagram}, which enables us to formulate a precise bulk-defect correspondence principle.

\begin{proposition} Let $\Gg$ be the canonical \'etale groupoid associated with the pattern $\Ll_0$ from Fig.~\ref{Fig:ConeNet}. Then $C^\ast(\Gg|_{\Ll_\infty}) \simeq C\big (\SM^1,C^\ast(\ZM^2) \big )$. Furthermore, $C^\ast(\Gg|_\Sigma) \simeq \KM(\ell^2(\Ll_0))$. 
\end{proposition}
\proof For any $x \in \Ll_{\infty}(\theta_0)$, we have $\Ll_{\infty}(\theta_0) - x = \Ll_{\infty}(\theta_0)$. Then the convolution~\eqref{Eq:Multi} reduces to the multiplication in the algebra $C^\ast(\ZM^2)$. The second statement is obvious.\qed

The algebra $C\big (\SM^1,C^\ast(\ZM^2) \big )$ is nuclear so the semi-group of extensions of this algebra by $\KM$ is actually a group \cite{KasparovJSM1981}. Then:

\begin{corollary} The bulk-defect correspondence  principle is encoded by the exact sequence
\begin{equation}\label{Eq:ExactSeq3}
\begin{diagram}
0 & \rTo & \KM & \rTo & C_r^\ast(\Gg) & \rTo & C\big (\SM^1,C^\ast(\ZM^2) \big ) & \rTo & 0,
\end{diagram}
\end{equation}
and, as such, by an element
\begin{equation}
[{\rm ext}]^1 \in KK^1 \Big (C\big (\SM^1,C^\ast(\ZM^2) \big ),\CM \Big ) \simeq \ZM^4.
\end{equation} 
\end{corollary}

We recall that the generators of 
$$KK^1 \Big (C\big (\SM^1,C^\ast(\ZM^2) \big ),\CM \Big ) \simeq KK^1 \big (C^\ast(\ZM^3),\CM \big )$$
are represented by a top Dirac operator $D_{\rm top} = \sum_{i=1}^3 \sigma_i \otimes X_i$, where $\sigma_i$'s are Pauli's matrices, and three lower Dirac operators, which involve only one direction and are in fact given by the components of the position operators $D_j = X_j$ \cite{ProdanRMP2016}.

\begin{proposition} The extension element $[{\rm ext}]^1$  contains the top generator $[D_{\rm top}]^1$.
\end{proposition}
\proof We will show numerically that $[u_{\rm top}]_1 \times [{\rm ext}]^1 \neq 0$. This can happen if and only if $[{\rm ext}]^1$ is a multiple of $[D_{\rm top}]^1$ (see \cite{ProdanSpringer2016}[Th.~5.7.1]).\qed

\begin{figure}[t]
\center
\includegraphics[width=\textwidth]{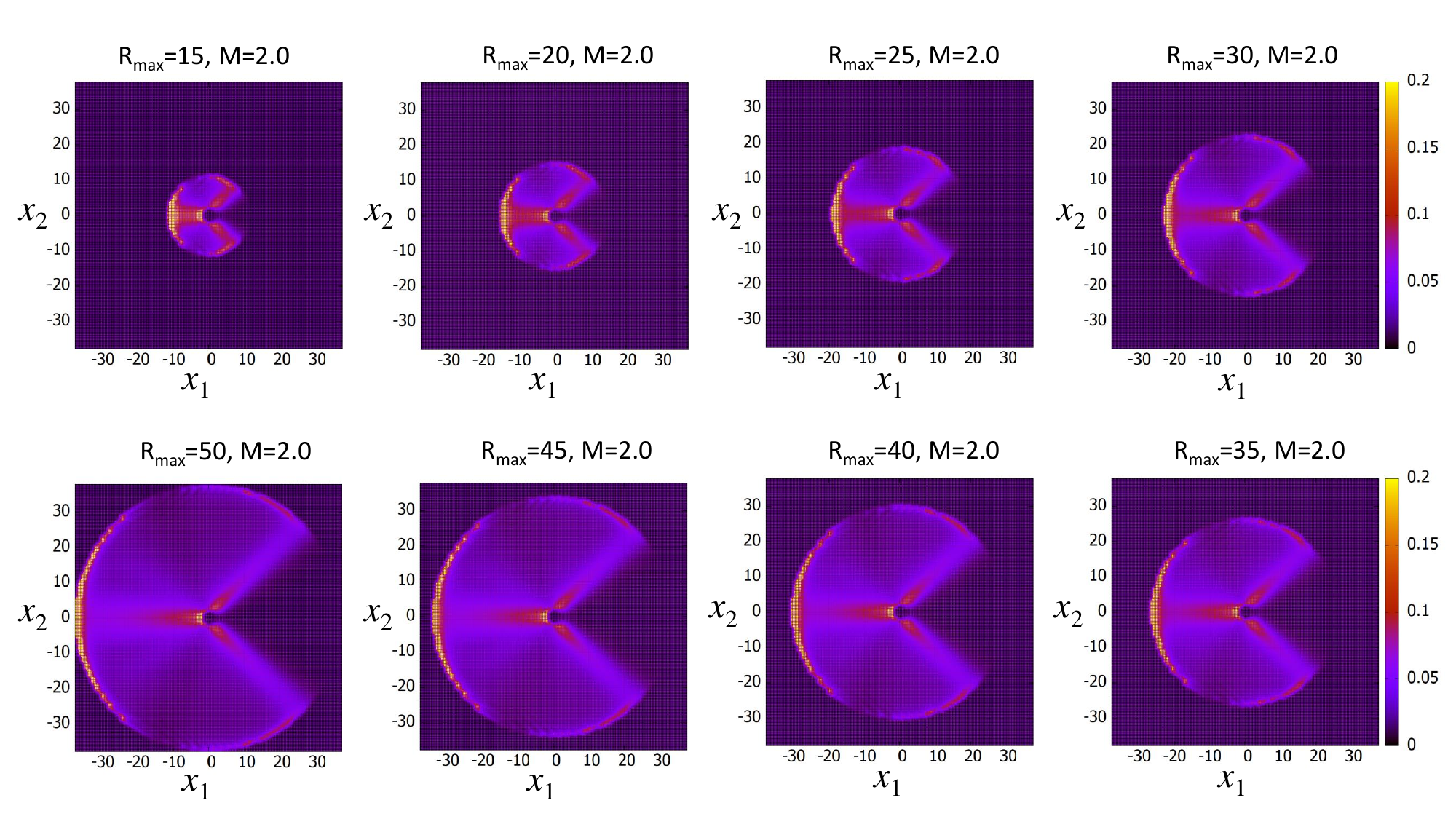}\\
  \caption{\small Local density of states at zero energy Eq.~\ref{Eq:LDOS}, rendered as function of $x_1$ and $x_2$ coordinates of a point $x \in \Ll_0$, $x = x_1 \hat e_x + x_2 \hat e_y + x_3 \hat e_z$.  
  }
 \label{Fig:ZeroModeVsSize}
\end{figure}

\section{A Numerical Demonstration}

The minimal model Hamiltonian on a square lattice which displays strong topological phases for class AIII ({\it i.e.} with chiral symmetry) in dimension 3 is defined on the Hilbert space $\CM^4 \otimes \ell^2(\ZM^3)$ and takes the form \cite{ProdanSpringer2016}[Sec.~2.3.3] 
\begin{equation}\label{Eq:Model1}
H_{\rm model} = \tfrac{1}{2\imath}\sum_{i=1}^3 \Gamma_i \otimes (S_i-S^\ast_i) + \Gamma_4 \otimes [M  + \tfrac{1}{2}\sum_{i=1}^3 (S_i + S^\ast_i)],
\end{equation}
where $S_i$'s are the shift operators of the 3-dimensional square lattice and $\Gamma$'s are $4 \times 4$ matrices supplying an irreducible representation of the Clifford algebra with four generators. The spectrum of the Hamiltonian displays a gap at the origin for all $M$ except for the critical values $M_c=\pm 3,\ \pm 1$. Furthermore, while varying $M$ continuously from $+\infty$ to $-\infty$, it is known that the $K_1$-class of the unitary operator $u_G$ associated to $H_{\rm model}$ takes the values $[1]_1$, $[u_{\rm top}]_1$,  $-2[u_{\rm top}]_1$, $[u_{\rm top}]_1$ and $[1]_1$, in this order. Here, $u_{\rm top}$ is the top generator of the $K_1$-group of $C^\ast(\ZM^3)$. The changes in the class of $u_G$ happens abruptly at the listed critical values $M_c$. We can transfer the model to the algebra $C\big (\SM^1,C^\ast(\ZM^2)\big)$ via the isomorphism $C^\ast(S_3) \simeq C(\SM^1)$, and we write the result explicitly, just for completeness ($\beta \in \SM^1$),
\begin{equation}\label{Eq:Model11}
\begin{aligned}
H_{\rm model}(\beta) & = \tfrac{1}{2\imath}\sum_{i=1,2} \Gamma_i \otimes (S_i-S^\ast_i) + \tfrac{1}{2} \Gamma_4 \otimes \sum_{i=1,2} (S_i + S^\ast_i) \\
& + \sin(\beta) \Gamma_3 \otimes 1 + \big (M+\cos(\beta)\big )\Gamma_4 \otimes 1.
\end{aligned}
\end{equation}

Our goal is to define a Hamiltonian on the pattern $\Ll_0$ shown in Fig.~\ref{Fig:ConeNet}, which converges to the model~\eqref{Eq:Model11} in the asymptotic regime $\Ll_0-x$, $|x| \rightarrow \infty$. We recall that one important assumption for the groupoid formalism is that the resonators are identical throughout the pattern and the connections between them should be entirely determined by the pattern. Hence, we will fix a representation for the Clifford matrices as
\begin{equation}
\Gamma_i =\begin{pmatrix} 0 & \sigma_i^\ast \\ \sigma_i & 0 \end{pmatrix}, \quad i =\overline{1,4},
\end{equation}
where $\sigma_4 =\imath I_{2 \times 2}$. We propose the following Hamiltonian:
\begin{equation}
\begin{aligned}
& H_{\Ll_0} = \sum_{x \in \Ll_0} \Big \{ \big(m-\tfrac{1}{2}\big)\Gamma_4 \otimes |x\rangle \langle x| + \tfrac{1}{2}\sum_{y \in \Ll_0} \chi_{xy}\Big [\imath \vec e_{yx} \cdot \vec \Gamma + \Gamma_4\Big] \otimes |y\rangle \langle x | \\
&  + \tfrac{1}{8}\sum_{z,y\in \Ll_0} \chi_{xyz}\Big [ \frac{(\hat e_y \cdot \vec e_{zyx})(\vec e_{zyx} \cdot \vec \Gamma)}{\cos(\xi)} + \frac{(\hat e_z \cdot \vec e_{zyx})(\hat e_x \cdot \vec e_{zyx})\Gamma_4}{\sin(\xi) \cos(\xi)} \Big ] \otimes |x\rangle \langle x| \Big \},
\end{aligned}
\end{equation} 
where $\vec \Gamma = \Gamma_1 \hat e_x + \Gamma_2 \hat e_y + \Gamma_3 \hat e_z$ and $\vec e_{yx} = y-x$, $\vec e_{zyx}=(z-x) \times (y-x)$. Also, $\chi_{xy}$ and $\chi_{xyz}$ are cut-off functions, which were chosen such that $x$, $y$ and $z$ are always nearest-neighbors. It is straightforward to check that the model converges to 
\begin{equation}\label{Eq:Model3}
\begin{aligned}
H_{\rm \Xi_\infty}(\theta_0) & = \tfrac{1}{2\imath}\sum_{i=1,2} \vec a_i(\theta_0) \cdot \vec \Gamma \otimes (S_i-S^\ast_i) + \tfrac{1}{2} \Gamma_4 \otimes \sum_{i=1,2} (S_i + S^\ast_i) \\
& + \sin(\theta_0/\alpha) \big (\vec a_1(\theta_0) \times \vec a_2(\theta_0) \big ) \cdot \vec \Gamma \otimes 1 + \big (M+\cos(\theta_0/\alpha)\big )\Gamma_4 \otimes 1,
\end{aligned}
\end{equation}
along the path defined by $\theta_0$ all the way to $\Xi_\infty$. Here, $S_i$ is the shift operator on $\Ll_\infty(\theta_0)$ in the $\hat a_i$ direction, $S_i|x\rangle = |x - \hat a_i\rangle$. As one can see, Eq.~\eqref{Eq:Model3} reproduces the topological model~\eqref{Eq:Model11}. According to the discussion at the beginning, the class of $u_G$ associated with the model~\eqref{Eq:Model3}, which belongs to $C^\ast(\Gg|_{\Xi_\infty}) \simeq C\big(\SM^1,C^\ast(\ZM^2)\big)$, becomes that of the top generator of the $K_1$-group exactly when $M$ becomes less than 3.

Fig.~\ref{Fig:RadLdosVsM} reports the local density of states ($\epsilon = 0.06$)
\begin{equation}\label{Eq:LDOS}
{\rm LDOS}(E,x) = {\rm Im} \, \langle x|(H_{\Ll_0}-E+ \imath \epsilon)^{-1}|x \rangle, \quad E \in \RM, \quad x \in \Ll_0,
\end{equation}
plotted as function of energy $E$ and $|x|$ for different values of $M$. The results reveal the emergence of zero-energy modes exactly at the critical value $M_c=3$. This confirms that the bulk-boundary correspondence is non-trivial, hence the extension module associated with the defect is indeed represented by a multiple of $[D_{\rm top}]^1$. The data reported in Fig.~\ref{Fig:RadLdosVsM} and additional un-reported data indicate that the zero-energy spectrum is not isolated.

Fig.~\ref{Fig:ZeroModeVsSize} reports the local density of states evaluated at $E=0$, rendered as function of $x_1$ and $x_2$ coordinates of $x \in \Ll_0$. In these simulations, the points with $|x| < 2.5$ have been removed from $\Ll_0$, just to demonstrate that the predicted topological zero modes are not affected. The simulations were performed for finite systems which include all $x \in \Ll_0$ with $|x| < R_{\rm max}$, with $R_{\max}$ specified above each panel of Fig.~\ref{Fig:ZeroModeVsSize}. The results supply strong evidence of a concentration of LDOS at $E=0$ near the core of the defect and that this concentration persists in the infinite size limit. According to the bulk-defect correspondence established in the previous section, this LDOS concentration cannot be removed by addition of any potential with support in a compact neighborhood of the tip of the cone. Let us specify that the pattern seen in Fig.~\ref{Fig:ZeroModeVsSize} can be change by re-parametrizing the circle $\Xi_\infty$.

\begin{figure}[t]
\center
\includegraphics[width=\textwidth]{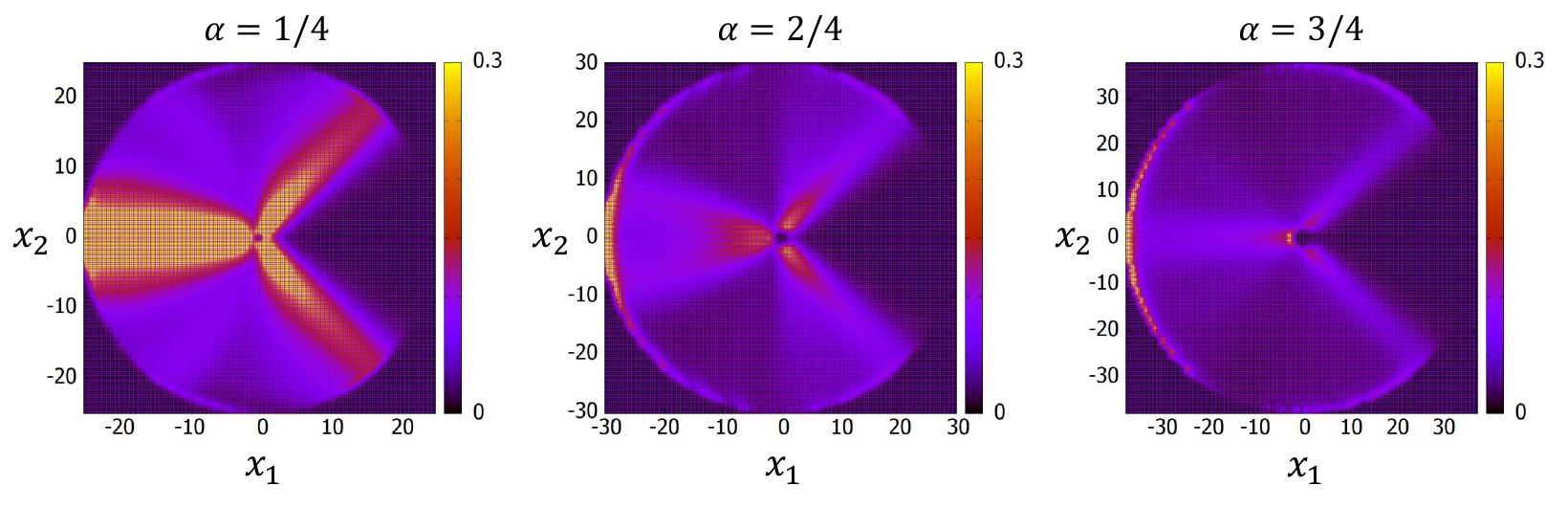}\\
  \caption{\small Mappings of the zero-energy LDOS for a sequence of defects. 
  }
 \label{Fig:DefectSequence}
\end{figure}

\section{Discussion}
\label{Sec:Discussion}

There is an entire sequence of topological defects for the 2-dimensional square lattice. They can be achieved by cutting the lattice along a seam and inserting or removing a quarter of the square lattice, leading to a deficit angle $n \frac{\pi}{2}$. There is no limit on how many quarters can be inserted, hence $n$ can take any negative integer number. However, we can only subtract at most three quarters, hence $n$ can only take the positive values 1, 2, 3. We have repeated the numerical simulations for these three cases, which correspond to $\alpha=3/4$, $2/4$ and $1/4$ in Eq.~\eqref{Eq:PattMap}, and the results are reported in Fig.~\ref{Fig:DefectSequence}. As one can see, the bulk-boundary correspondence persists and, in fact, the zero-energy LDOS becomes more intense with $n$.

One can define a natural addition of two defects, by opening a seam for each of the defects and re-gluing appropriately. This correlates well with the semi-group structure of extensions, particularly, with addition of Busby invariants \cite{JensenBook}[p.~98]. It is very desirable to make this connection explicit. Furthermore, a pair $(n,-n)$ of defects can be nucleated from a pristine lattice without affecting $\Xi_\infty$, which is just a point for the pristine lattice. Hence, $n$ and $-n$ defects can be thought as inverse to each other. The extension semi-group is a group in our setting and such pairs of mutually inverse extensions exist for any $n$. However, to go beyond $n=3$, one needs to use the stabilization of the algebras and it is not clear to us if such extensions can be represented as a groupoid algebra extensions of the type encountered here. Nevertheless, we conjecture that the $n$ defect is associated with the element $n[D_{\rm top}]^1$ of the Kasparov group $KK^1(C^\ast(\ZM^3),\CM)$.

Let us end with the note that, so far, the discussion involved only the additive structure of $KK$-theory. It will be interesting to understand if the Kasparov product plays any role here. Perhaps, the spaces of units that look as in Fig.~\ref{Fig:GenTrans} can be split in ways that lead to Kasparov elements that can be multiplied. This could be interpreted as fusing the defects, hence supplying examples where a full theory of topological defects can be formalized inside Kasparov's formalism.

\end{document}